\documentclass[letterpaper, 10 pt, journal, twoside]{IEEEtran}

\title{\LARGE \bf
Piecewise semi-ellipsoidal control invariant sets
}

\author{Beno\^it Legat and Sa\v{s}a V. Rakovi\'c and Rapha\"el M. Jungers
\thanks{B. Legat is a FNRS Research Fellow.
Address: ICTEAM institute, UCLouvain, Louvain-la-Neuve, Belgium.
{\tt\small benoit.legat@uclouvain.be}}
\thanks{S. V. Rakovi\'c is with Beijing Institute of Technology, Beijing, China. {\tt\small sasa.v.rakovic@gmail.com}}%
\thanks{R. M. Jungers is a FNRS honorary Research Associate. He is supported by the Walloon Region and the Innoviris Foundation.
Address: ICTEAM institute, UCLouvain, Louvain-la-Neuve, Belgium.
{\tt\small raphael.jungers@uclouvain.be}}%
}

\usepackage{comment}

\usepackage{caption}
\usepackage{subcaption}

\usepackage{custom}

\usepackage{mythm}

\usepackage{mathtools}
\newcommand{\verteq}{\rotatebox{90}{$=$}}

\renewcommand{\SymK}{\mathbb{S}^n}
\renewcommand{\Psd}{\mathbb{S}_+^n}
\renewcommand{\polar}[1]{{#1}^*}

\DeclareMathOperator{\gauge}{g}
\DeclareMathOperator{\suppfun}{h}

\begin{document}

\maketitle
\thispagestyle{empty}
\pagestyle{empty}

\begin{abstract}

  Computing control invariant sets is paramount in many applications.
  The families of sets commonly used for computations are ellipsoids and polyhedra.
  However, searching for a control invariant set over the family of ellipsoids is
  conservative for systems more complex than unconstrained linear time invariant systems.
  Moreover, even if the control invariant set may be approximated arbitrarily
  closely by polyhedra, the complexity of the polyhedra may grow rapidly
  in certain directions.
  An attractive generalization of these two families are piecewise semi-ellipsoids.
  We provide in this paper a convex programming approach for computing control
  invariant sets of this family.

\end{abstract}

\section{Introduction}

\IEEEPARstart{C}omputing control invariant sets is paramount in many applications~\cite{blanchini2015set}.
The existence of a nontrivial control invariant set of a linear time-invariant (LTI) control system is equivalent
to the stability of its uncontrollable subspace (which is readily accessible
in its Controllability Form) \cite[Section~2.4]{wonham1974linear}.
Indeed, the eigenvalues of its controllable subspace can be fixed to any value
by a proper choice of linear state feedback.
The resulting autonomous system is stable hence an invariant ellipsoid can be
determined by solving a system of linear equations~\cite{Liapounoff1907}.
This set is also control invariant for the control system.

While searching for an ellipsoidal invariant set is not conservative for autonomous LTI dynamics,
it is no longer the case for \emph{uncertain} or \emph{switched} systems \cite{petersen1985quadratic}.
Furthermore, it is often desirable to find a control invariant set of maximal
volume (or which is maximal in some direction~\cite{ahmadi2018robust}).
For such problem, the method of fixing the eigenvalues of the controllable subspace detailed above is not suitable as it does not
take any volume consideration into account.
More importantly, the maximal
invariant set may not be an ellipsoid and may not be rendered invariant via
a linear control.
In this paper, we consider both a switched dynamics and polyhedral constraints on the state space.

The maximal control invariant set of a linear control system
can be obtained as the limit of a fixed point iteration,
since it is the maximal fixed point of the standard viability kernel algorithm.
When the set of state constraints and the set of input constraints are polyhedra,
each iterate is a polyhedron but the limit may not be polyhedral.
Computing the next iterate requires computing the projection of a preimage of the previous iterate
(see e.g., the procedure p.~201 in \cite{blanchini2015set}).
While computing a preimage of a polyhedron given its \emph{H-representation} (see \eqref{eq:hrep})
is computationally cheap and is represented by the same number of halfspaces,
its projection may increase the complexity of the representation by generating numerous new halfspaces,
e.g. by Fourier-Motzkin elimination.

In \secref{piece}, we define the family of piecewise semi-ellipsoids and
review the literature for computing invariant sets of this family for autonomous systems.
In \secref{alg}, we review the common algebraic approach for computing invariant sets
and discuss the challenge of computing piecewise semi-ellipsoidal control invariant sets for control systems with this approach.
In \secref{geo}, we detail a geometric approach for computing
control invariant sets and show that this allows to compute piecewise semi-ellipsoidal control invariant sets.
We present the approach for discrete-time control linear systems for simplicity and then
show how to generalize it for switched systems.
We illustrate the technique on a simple example.

\paragraph*{Reproducibility}
The code used to obtain the results is published on codeocean~\cite{legat2020piecewiseCO}.
The algorithms are part of the SetProg package~\cite{legat2019set} in Julia~\cite{bezanson2017julia}
which solves a set program for given templates.
The semidefinite programs formulated by SetProg for the ellipsoidal and piecewise semi-ellipsoidal templates are solved by Mosek v8~\cite{mosek2017mosek81034} through JuMP~\cite{dunning2017jump}.
The polyhedral computations needed for the fixed point iterations of the standard viability kernel algorithm are carried out using the Polyhedra package~\cite{legat2020polyhedra}.
The linear programs formulated by Polyhedra to remove redundant elements of a H-representation
are solved by Mosek v8~\cite{mosek2017mosek81034} through MathOptInterface~\cite{legat2020mathoptinterface}.

\subsection{Notation}

The set of nonnegative real number is denoted $\R_+$,
the set of symmetric $n \times n$ square matrices is denoted $\SymK$
and the set of positive semidefinite matrices of $\SymK$ is denoted $\Psd$.
The notation $A \succeq B$ denotes that $A - B \in \Psd$.
The notation $A^\pseudoinv$ denotes the pseudoinverse of the matrix $A \in \Psd$.

Given a matrix $A \in \R^{n_1 \times n_2}$ and a set $\mathcal{S} \subseteq \R^{n_2}$,
the \emph{image} of $\mathcal{S}$ under $A$ is
\(
  A\mathcal{S} = \{\, Ax \mid x \in \mathcal{S} \,\}.
\)
For a scalar number $\gamma \in \mathbb{R}$, $\gamma \mathcal{S}$ denotes $\gamma I_{n_2} \mathcal{S}$
where $I_{n_2}$ is the identity matrix of dimension $n_2$.

Given a matrix $A \in \R^{n_1 \times n_2}$ and a set $\mathcal{S} \subseteq \R^{n_1}$,
the \emph{preimage} of $\mathcal{S}$ under $A$ is
\(
  A^{-1}\mathcal{S} = \{\, x \mid Ax \in \mathcal{S} \,\}
\)
Note that $A$ does not need to be invertible in this definitions.
We also denote the preimage of $A^\Tr$ as $A^{-\Tr}$.

Given a linear subspace $\mathcal{V}$, $\mathcal{V}^\perp$ denotes its orthogonal
subspace and $P_{\mathcal{V}}$ denotes the projection onto $\mathcal{V}$.

Hyperplanes and halfspaces are respectively denoted by
\(
  \hhyperplane{a} = \{\, x \mid a^\Tr x = 0 \,\}
\)
and
\(
  \hhalfspace{a} = \{\, x \mid a^\Tr x \ge 0 \,\}
\).
A \emph{H-representation} of a polyhedral cone $\mathcal{P} \subseteq \R^n$ 
is a set of vectors $a_i \in \R^n$ such that
\begin{equation}
  \label{eq:hrep}
  \mathcal{P} = \{\, x \mid a_i^\Tr x \ge 0, \forall i \,\} = \bigcap_i \hhalfspace{a_i}.
\end{equation}
In other words, it is the representation of a polyhedron as a finite intersection of halfspaces.
The \emph{affine hull} of a polyhedron $\mathcal{P}$, denoted by $\aff{\mathcal{P}}$, is the smallest affine space
containing $\mathcal{P}$.
The \emph{dimension} of $\mathcal{P}$, denoted $\dim(\mathcal{P})$ is the dimension of its affine hull.
The convex hull of the union of convex sets $(\mathcal{S})_{i \in I}$ is denoted $\conv_{i \in I} \mathcal{S}_i$.
The \emph{V-representation} of a polytope $\mathcal{P} \subseteq \R^n$, denoted $\Vrep(\mathcal{P})$,
is the minimal set of points $v_i$ such that $\mathcal{P} = \conv_i \{v_i\}$.

\section{Piecewise semi-ellipsoids}
\label{sec:piece}

The conservatism of ellipoidal control invariant sets and the complexity of the representation
of polyhedral control invariant sets has motivated the search for alternative templates of sets.
The template we study in this paper is the family of piecewise semi-ellipsoids.
These sets may either be defined as the 1-sublevel sets of piecewise quadratic forms
or as sets with a piecewise semi-Euclidean \emph{Minkowski function}.
\begin{mydef}[Minkowski function]
  \label{def:gauge}
  We define the \emph{gauge} or \emph{Minkowski function} of a closed convex set $\mathcal{S}$ containing the origin as:
  \( \gauge(\mathcal{S}, x) = \min_\gamma \{\,\gamma : x \in \gamma \mathcal{S}, \gamma \ge 0\,\}\), for $x \in \R^n$.
\end{mydef}
\begin{myrem}
  \label{rem:gaugeinf}
  As we do not assume boundedness of the set in \defref{gauge}, the Minkowski function can be zero for nonzeros vectors $x$.
  Moreover, as we do not assume that the set has nonempty interior and that the origin is in the interior of the set,
  the Minkowski function may be infinite for some vectors $x$.
  The image set of Minkowski functions is the set of extended reals for these reasons.
  As we compare Minkowski functions in this paper, we adopt the
  convention $\infty = \infty$ so that the inequality $\gauge(\mathcal{S}_1, x) \ge \gauge(\mathcal{S}_2, x)$ holds when the value on both sides is $\infty$.
  See \cite[Example~3.50]{rockafellarvariational} for more details.
\end{myrem}
The set $\mathcal{S}$ is the 1-sublevel set of its Minkowski function and the Minkowski function
is the only sublinear function that has this property for the set $\mathcal{S}$.
For instance, the ellipsoids are the 1-sublevel sets of quadratic forms but quadratic forms are not sublinear.
The Minkowski functions of ellipsoids are the square root of quadratic forms.
\begin{mydef}[Conic partition]
  A \emph{conic partition} of $\R^n$ is a set of $m$ polyhedral cones $(\mathcal{P}_i)_{i=1}^m$
  with nonempty interior
  such that for all $i \neq j$, $\dim(\mathcal{P}_i \cap \mathcal{P}_j) < n$
  and $\cup_{i=1}^m \mathcal{P}_i = \R^n$.
\end{mydef}
Given a conic partition, we use the notation
\(
  \mathcal{N} = \{\, (i, j) \mid \dim(\mathcal{P}_i \cap \mathcal{P}_j) = n - 1 \,\}.
\)
For each $(i, j) \in \mathcal{N}$,
the affine hull, of $\mathcal{P}_i \cap \mathcal{P}_j$ is a hyperplane.
We denote by $n_{ij}$ the normal of this hyperplane directed towards $\mathcal{P}_i$,
i.e., such that $\mathcal{P}_i \subseteq \hhalfspace{n_{ij}}$.
\begin{mydef}[Piecewise semi-ellipsoids]
  \label{def:piece}
  A closed convex set $\mathcal{S} \subseteq \R^n$ containing the origin is said to be \emph{piecewise semi-ellipsoidal} if there exists
  a conic partition $(\mathcal{P}_i)_{i=1}^m$ of $\R^n$
  and symmetric positive semidefinite matrices $Q_i \in \Psd$ for $i = 1, \ldots, m$,
  such that
  \begin{align}
    \notag
    \gauge(\mathcal{S}, x) & =
      \sqrt{x^\Tr Q_i x} \quad \text{ if } x \in \mathcal{P}_i,\\
    \label{eq:piececont}
    x^\Tr Q_i x & = x^\Tr Q_j x, \quad \forall i, j \in \mathcal{N}, x \in \mathcal{P}_i \cap \mathcal{P}_j,\\
    \label{eq:piececonv}
    n_{ij}^\Tr Q_i x & \ge n_{ij}^\Tr Q_j x, \quad \forall i, j \in \mathcal{N}, x \in \mathcal{P}_i \cap \mathcal{P}_j. 
  \end{align}
\end{mydef}
Note that without \eqref{eq:piececont} and \eqref{eq:piececonv}, the 1-sublevel set of
the piecewise semi-Euclidean function may be non-convex.
The condition \eqref{eq:piececont} ensures the continuity of the function
and \eqref{eq:piececonv} ensures its convexity.
Indeed, \eqref{eq:piececont} and \eqref{eq:piececonv} ensures that $\gauge(\mathcal{S}, x) = \max(\sqrt{x^\Tr Q_i x}, \sqrt{x^\Tr Q_i x})$ in the neighbourhood $\mathcal{P}_i \cap \mathcal{P}_j$ contained in the interior of $\mathcal{P}_i \cup \mathcal{P}_j$.
By continuity of the square root of quadratic forms $\sqrt{x^\Tr Q_i x}$, $\gauge(\mathcal{S}, x) = \max_{i \in I} \sqrt{x^\Tr Q_i x}$ in the neighbourhood of $\cap_{i \in I} \mathcal{P}_i$ contained in the interior of $\cup_{i \in I} \mathcal{P}_i$ for any subset $I \subseteq \{1,\ldots,n\}$.

This family generalizes both ellipsoids and polyhedra.
Indeed, if $m = 1$ and $\mathcal{P}_1 = \R^n$, we recover the family of ellipsoids and
if the matrices $Q_i$ are rank-1, we recover the family of polyhedra.

Given a piecewise semi-ellipsoid, its polar is also piecewise semi-ellipsoidal
but the conic partition of the polar depends on the matrices $Q_i$.
Indeed, the Minkowski function of the polar is given by
$\gauge(\polar{\mathcal{S}}, x) = \sqrt{2(\gauge(\mathcal{S}, x)^2 / 2)^*}$ (see \cite[Corollary~15.3.2]{rockafellar2015convex})
which is piecewise semi-ellipsoidal as $\gauge(\mathcal{S}, x)^2 / 2$ is piecewise quadratic
and the conjugate of a piecewise quadratic function is a piecewise quadratic function; see \cite[Theorem~11.14]{rockafellarvariational}.
The Minkowski function of the polar set has the closed form expression given by \propref{polarpiece}.

\begin{myprop}
  \label{prop:polarpiece}
  Given a piecewise semi-ellipsoid $\mathcal{S}$,
  as defined in \defref{piece}, the polar set $\polar{\mathcal{S}}$
  is the piecewise semi-ellipsoid represented by
  the conic partition made of the polyhedral cones $Q_i\mathcal{P}_i$
  with matrices $Q_i^\pseudoinv$ for $i = 1, \ldots, m$
  and the polyhedral cones
  \begin{equation}
    \label{eq:newpart}
    \conv_{i \in I} Q_i \bigcap_{i \in I} \mathcal{P}_i
  \end{equation}
  with matrices
  \(
    E_I^\Tr(E_I Q_i E_I^\Tr)^\pseudoinv E_I
  \)
  where\footnote{We have $E_I Q_i E_I^\Tr = E_I Q_j E_I^\Tr$ for any $i, j \in I$ by \eqref{eq:piececont} so the matrix is independent on the $i \in I$ chosen.} $i \in I$ and $E_I = P_{\aff{\cap_{i \in I} \mathcal{P}_i}}$
  for any subset $I$ of $\{1, \ldots, m\}$.
  \begin{proof}
    Given a point $x$ in the intersection of the boundary of $\mathcal{S}$ and the interior of $\mathcal{P}_i$, the tangent cone to $\mathcal{S}$ at $x$
    contains a single direction given by $Q_i x$.
    Therefore, the value of the support function is given by $y^\Tr Q_i^\pseudoinv y$ for
    $y \in Q_i\mathcal{P}_i$. 
    Given a subset $I$ of $\{1, \ldots, m\}$, by \eqref{eq:piececont},
    the matrix $E_I Q_i E_I^\Tr$ is identical for all $i \in I$.
    For a point $x$ in the intersection of the boundary of $\mathcal{S}$ and $\cap_{i \in I} \mathcal{P}_i$,
    by \eqref{eq:piececonv}, the tangent cone to $\mathcal{S}$ at $x$
    is the conic hull of the vectors $Q_i x$ for each $i \in I$.
    Therefore, the value of the support function is given by $y^\Tr E_I^\Tr(E_I Q_i E_I^\Tr)^\pseudoinv E_I y$ for $y$ in this cone.
  \end{proof}
\end{myprop}

\begin{myrem}
  The conic partition for the polar set created by \propref{polarpiece}
  seems to contain many polyhedral cones.
  This seems surprising as the polar operation is an involution for closed convex sets containing the origin.
  In fact, many polyhedral cones of the conic partition created can be dropped without
  changing the set as they are not full-dimensional.
  For instance, the pieces of the partition created in \eqref{eq:newpart}
  are only full-dimensional in case the subdifferential of $\gauge(\mathcal{S}, x)$ is not a singleton for all $x \in \cap_{i \in I} \mathcal{P}_i$.
  That is, if \eqref{eq:piececonv} is not satisfied with equality for each pair of $i, j \in I$.
\end{myrem}

\exemref{piecequad} illustrates the computation of the polar of a piecewise semi-ellipsoid
with \propref{polarpiece}.

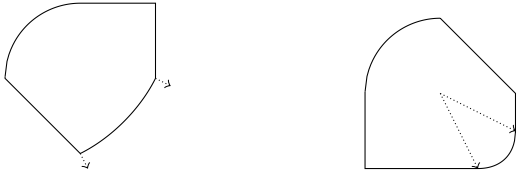
\begin{figure}[!ht]
  \centering
  \begin{subfigure}[b]{0.23\textwidth}
    \centering
    \begin{tikzpicture}
      \draw 
      (-1, 0.02) to (-1, 0) to (0, -1)
      plot[domain=0:1, samples=40] ({\x}, {\x/2 - sqrt(4 - 3 * \x^2) / 2} )
      (1, 0) to (1, 1) to (0, 1)
      plot[domain=0:1, samples=40] ({-\x}, {sqrt(1 - \x^2)} );
      \draw [densely dotted, ->] (0, -1) to (0.1, -1.2);
      \draw [densely dotted, ->] (1, 0) to (1.2, -0.1);
    \end{tikzpicture}
    \caption{Set $\mathcal{S}$ whose Minkowski function is defined by \eqref{eq:pieceprimal}.}
    \label{fig:primal}
  \end{subfigure}%
  \hfill%
  \begin{subfigure}[b]{0.23\textwidth}
    \centering
    \begin{tikzpicture}
      \draw 
      (-1, 0.02) to (-1, -1) to (0.5, -1)
      plot[domain=0.5:1, samples=40] ({\x}, {-\x/2 - sqrt(3 - 3 * \x^2) / 2} )
      (1, -0.52) to (1, 0) to (0, 1)
      plot[domain=0:1, samples=40] ({-\x}, {sqrt(1 - \x^2)} );
      \draw [densely dotted, ->] (0, 0) to (0.5, -1);
      \draw [densely dotted, ->] (0, 0) to (1, -0.5);
    \end{tikzpicture}
    \caption{Set $\polar{\mathcal{S}}$ whose Minkowski function is defined by \eqref{eq:piecepolar}.}
    \label{fig:polar}
  \end{subfigure}
  \caption{Illustration for sets $\mathcal{S}$ and $\polar{\mathcal{S}}$ defined in \exemref{piecequad}.}
  \label{fig:example}
\end{figure}

\begin{myexem}
  \label{exem:piecequad}
  Consider the piecewise semi-ellipsoid whose Minkowski function is defined by
  \begin{equation}
    \label{eq:pieceprimal}
    \gauge(\mathcal{S}, x) =
    \begin{cases}
      |x_1| & \text{ if } 0 \le x_2 \le x_1,\\
      |x_2| & \text{ if } 0 \le x_1 \le x_2,\\
      \sqrt{x_1^2 + x_2^2} & \text{ if } x_1 \le 0 \le x_2,\\
      |x_1 + x_2| & \text{ if } x_1, x_2 \le 0,\\
      \sqrt{x^\Tr Q_5 x} & \text{ if } x_2 \le 0 \le x_1.
    \end{cases}
  \end{equation}
  where $x^\Tr Q_5 x = x_1^2 - x_1x_2 + x_2^2$.
  The polar set is also piecewise semi-ellipsoidal
  and its Minkowski function is given by:
  \begin{equation}
    \label{eq:piecepolar}
    \gauge(\polar{\mathcal{S}}, x) =
    \begin{cases}
      |x_1 + x_2| & \text{ if } 0 \le x_1, x_2,\\
      \sqrt{x_1^2 + x_2^2} & \text{ if } x_1 \le 0 \le x_2,\\
      |x_1| & \text{ if } x_1 \le x_2 \le 0,\\
      |x_2| & \text{ if } x_2 \le x_1, 2x_1 + x_2 \le 0,\\
      \sqrt{x^\Tr Q_5^{-1} x} & \text{ if } 2x_1 + x_2 \ge 0, x_1 + 2x_2 \le 0,\\
      |x_1| & \text{ if } x_1 + 2x_2 \ge 0, x_2 \ge 0.
    \end{cases}
  \end{equation}
  where $x^\Tr Q_5^{-1} x = (4/3) \cdot (x_1^2 + x_1x_2 + x_2^2)$.

  Note that the partition $\mathcal{P}_i$ of $\polar{\mathcal{S}}$ does
  not only depend on the conic partition of $\mathcal{S}$, it also depends on the
  value of the matrices $Q_i$.

  For instance, the cone defined by $2x_1 + x_2 \ge 0, x_1 + 2x_2 \ge 0$
  is the conic hull of the gradient of $\sqrt{x_1^2 - x_1x_2 + x_2^2}$ evaluated
  at $(0, -1)$ and $(1, 0)$ as shown in dotted arrows in \figref{example}.
  The gradients are obtained by multiplying $(0, -1)$ and $(1, 0)$
  by the matrix $Q_5$.
  This illustrates why the partition of the polar is the image of
  the original partition under $Q_5$.
\end{myexem}


The study of piecewise quadratic Lyapunov functions started with continous-time autonomous piecewise affine systems in~\cite{johansson1998computation}.
The authors set the polyhedra $\mathcal{P}_i$ to the pieces used to define the system.
This is generalized to piecewise polynomial Lyapunov functions for continuous-time autonomous piecewise polynomial systems in~\cite{prajna2003analysis}.
The authors mention in \cite[Section~3.2.1]{prajna2003analysis} that refining the partition is not obvious for non-planar systems.
They suggest as an alternative to increase the degree of the polynomials.

\section{Algebraic approach}
\label{sec:alg}

We argue in this section that there is a fundamental challenge in computing piecewise quadratic Lyapunov functions for control system.
First, we review the computation of invariant ellipsoids for autonomous systems.

\subsection{Autonomous systems}

An ellipsoid $\mathcal{S}$, i.e. with Minkowski function
\( \gauge(\mathcal{S}, x) = \sqrt{x^\Tr Q x} \), is invariant for a discrete-time linear autonomous system
\begin{equation}
  \label{eq:autonsys}
  x_{k+1} = Ax_k
\end{equation}
if for all $x \in \R^n$ such that $x^\Tr Q x \leq 1$, we have $x^\Tr A^\Tr Q A x \leq 1$.
This condition is equivalent to the following Linear Matrix Inequality (LMI):
\begin{equation}
  \label{eq:uncontrol_LMI}
  Q \succeq A^\Tr Q A.
\end{equation}

A piecewise semi-ellipsoid is invariant for the same system if
\begin{quote}
  for all $i, j = 1, \ldots, m$, for all $x \in \mathcal{P}_i \cap A^{-1}\mathcal{P}_j$,
  we have $x^\Tr (Q_i - A^\Tr Q_j A) x \geq 0$.
\end{quote}
When $\mathcal{P}_i \cap A^{-1}\mathcal{P}_j$ is the nonnegative orthant, this is equivalent to the copositivity of $Q_i - A^\Tr Q_j A$.
While checking copositivity of a matrix is co-NP-complete~\cite{murty1987some},
sufficient conditions such as \propref{copos} can be encoded as a LMI.

\begin{myprop}[{\cite[Section~3.6.1]{blekherman2012semidefinite}}]
  \label{prop:copos}
  Consider a polyhedron $\mathcal{P}$ with H-representation $(a_i)_{i=1}^k$
  and a symmetric matrix $Q \in \SymK$.
  If there exists nonnegative $\lambda_{ij} \in \R_+$ such that
  \( Q - \sum_{i \neq j} \lambda_{ij} (a_ia_j^\Tr + a_ja_i^\Tr) \)
  is positive semidefinite then for all $x \in \mathcal{P}$, we have
  $x^\Tr Q x \ge 0$.
\end{myprop}

\begin{myrem}
  \label{rem:copos}
  While \propref{copos} only provides a sufficient condition,
  a necessary condition can be obtained using a hierarchy of
  semidefinite programs of increasingly larger size~\cite[Chapter~5]{parrilo2000structured}.
\end{myrem}

\subsection{Control systems}

We consider a linear control system with state constraints but no
input constraints\footnote{Note that this is without loss of generality as if there are input constraints, we can consider a lifted system where the input constraints are moved to state constraints; see \cite[Section~2.2]{legat2020sum} for more details.}
\begin{equation}
  \label{eq:controlsys}
  x_{k+1} = Ax_k + Bu_k, x_k \in \mathcal{X}.
\end{equation}
We say that a set $\mathcal{S} \subseteq \mathcal{X}$ is \emph{control invariant} for the system if
\begin{equation}
  \label{eq:controlinv}
  \forall x \in \mathcal{S}, \exists u : Ax + Bu \in \mathcal{S}
\end{equation}

If the Minkowski function is \( \gauge(\mathcal{S}, x) = \sqrt{x^\Tr Q x} \),
i.e. $\mathcal{S}$ is an ellipsoid,
then \eqref{eq:controlinv} is rewritten into
\begin{equation}
  \label{eq:existential}
  x^\Tr Q x \leq 1 \Rightarrow \exists u, (A x + Bu)^\Tr Q (A x + Bu) \leq 1.
\end{equation}
The control term $u$, or more precisely the existential quantifier $\exists$ prevents the condition to be rewritten into an LMI directly.

There is a well known technique to circumvent the presence of the existential
quantifier $\exists$ in \eqref{eq:existential}, which allows to formulate the search for an ellipsoidal
control invariant set of control linear discrete systems as a semidefinite
program.
This is described in details in \cite[Section~7.2.2]{boyd1994linear} for continuous-time
and in \cite[Section~4.4.2]{blanchini2015set} or \cite[Section~2.2.1]{blekherman2012semidefinite} for discrete-time.
We describe this technique in the following paragraph to highlight the challenge to generalize it for piecewise semi-ellipsoids.

Fixing the control to a linear state feedback $u(x) = Kx$ for some matrix $K$ allows
to use the condition \eqref{eq:uncontrol_LMI} for the linear autonomous system $x_{k+1} = (A + BK) x_k$.
Using the Schur lemma, the invariance condition can be formulated
as a \emph{Bilinear Matrix Inequality} (BMI) which is NP-hard to solve in general \cite{toker1995np}.
While the matrix inequality is bilinear in $K$ and $Q$,
a clever algebraic manipulation allows to reformulate it as a
LMI in the decision variables $Q^{-1}$ and $Y := KQ^{-1}$.
The LMI is
\begin{equation}
  \label{eq:control_LMI}
  \begin{bmatrix}
    Q^{-1} & Q^{-1}A^\Tr + Y^\Tr B^\Tr\\
    AQ^{-1} + BY & Q^{-1}
  \end{bmatrix} \succeq 0.
\end{equation}

For piecewise semi-ellipsoids, given a state vector $x \in \mathcal{P}_i$,
the next iterate is $(A + BK)x$.
Therefore, we need to somehow use \eqref{eq:control_LMI} with $Q_i$ and $Q_j$
on the polyhedra $\mathcal{P}_i \cap (A + BK)^{-1}\mathcal{P}_j$.
However, because of the reformulation into the decision variable $Y$,
$K$ is not a decision variable of the semidefinite program.
It is therefore unclear how to formulate the control invariance of a piecewise semi-ellipsoid
for a linear control system.

\section{Geometric approach}
\label{sec:geo}

We first introduce \emph{support functions}.

\begin{mydef}[Support function]
  We define the \emph{support function} of a nonempty closed convex set $\mathcal{S} \subseteq \R^n$ as:
  \( \suppfun(\mathcal{S}, y) = \sup_{x \in \mathcal{S}} y^\Tr x\), for $y \in \R^n$.
\end{mydef}
\begin{myrem}
  \label{rem:supportinf}
  Similarly to \remref{gaugeinf}, the image of support functions is the extended reals
  and we adopt the same convention as \remref{gaugeinf} for comparing infinite values.
\end{myrem}

As shown by the following proposition, the support function is the Minkowski function of the polar set.

\begin{myprop}[{\cite[Theorem~14.5]{rockafellar2015convex}
  }]
  \label{prop:gaugesupp}
  Given a closed convex set $\mathcal{S} \subseteq \R^n$ containing the origin,
  for all $x\in \mathbb{R}^n$, we have
  \(\gauge(\mathcal{S}, x) = \suppfun(\polar{\mathcal{S}}, x).\)
\end{myprop}

We have the following property for the Minkowski functions of linear preimages and support functions of linear images.
\begin{myprop}
  \label{prop:gA}
  Given a matrix $A \in \R^{n_1 \times n_2}$ and a closed convex set $\mathcal{S} \subseteq \R^{n_1}$ containing the origin,
  for all $x \in \R^{n_2}$,
  the following holds:
  \begin{equation}
    \label{eq:gA}
    \gauge(A^{-1}\mathcal{S}, x) = \gauge(\mathcal{S}, Ax)
  \end{equation}
\end{myprop}
\begin{myprop}[{\cite[Corollary~11.24(c)]{rockafellarvariational} or \cite[Corollary~16.3.1]{rockafellar2015convex}}]
  \label{prop:hA}
  Given a matrix $A \in \R^{n_1 \times n_2}$ and a nonempty closed convex
  set $\mathcal{S} \subseteq \R^{n_2}$, for all $y \in \R^{n_1}$,
  the following holds:
  \begin{equation}
    \label{eq:hA}
    \suppfun(A\mathcal{S}, y) = \suppfun(\mathcal{S}, A^\Tr y)
  \end{equation}
\end{myprop}

The following properties show the relation between the linear image of a convex set and its polar.
\begin{myprop}[{\cite[Corollary~16.3.2]{rockafellar2015convex}}]
  \label{prop:podu}
  For any convex set $\mathcal{S}$ and linear map $A$,
  we have
  \(
    \polar{(A\Csetvar)} = A^{-\Tr} \polar{\Csetvar}.
  \)
\end{myprop}

\subsection{Autonomous systems}

The invariance condition of a set $\mathcal{S}$ for the linear autonomous system \eqref{eq:autonsys} can be written as
\begin{equation}
  \label{eq:ASS}
  A\mathcal{S} \subseteq \mathcal{S}.
\end{equation}
We see with \propref{image} that this is equivalent to
\begin{equation}
  \label{eq:SAS}
  \mathcal{S} \subseteq A^{-1}\mathcal{S}.
\end{equation}

\begin{myprop}
  \label{prop:image}
  Consider a matrix $A \in \R^{n_1 \times n_2}$, a set $\mathcal{S} \subseteq \R^{n_2}$ and a set $\mathcal{T} \subseteq \R^{n_1}$.
  The inclusion $A\mathcal{S} \subseteq \mathcal{T}$ holds if and only if
  the inclusion $\mathcal{S} \subseteq A^{-1}\mathcal{T}$ holds.
  \begin{proof}
    The inclusion $A\mathcal{S} \subseteq \mathcal{T}$
    is equivalent to $x \in \mathcal{S} \Rightarrow Ax \in \mathcal{T}$.
    Since $Ax \in \mathcal{T}$ is equivalent to $x \in A^{-1}\mathcal{T}$,
    we have the desired result.
  \end{proof}.
\end{myprop}

\begin{myrem}
  \label{rem:image}
  Note that $\mathcal{S} \subseteq A\mathcal{T}$
  does not imply $A^{-1}\mathcal{S} \subseteq \mathcal{T}$.
  For instance, if
  $A = \begin{bmatrix} 1 & 0
  \end{bmatrix}$,
  $\mathcal{S} = [-1, 1]$ and $\mathcal{T} = [-1, 1]^2$ then
  $\mathcal{S} = [-1, 1] = A\mathcal{T}$ but
  $A^{-1}\mathcal{S} = [-1, 1] \times \R \not\subseteq \mathcal{T}$.
\end{myrem}

We can formulate invariance inequalities in terms of either the support (resp. gauge) function of $\mathcal{S}$
with \eqref{eq:suppASS} (resp. \eqref{eq:gaugeSAS}) or
we can write them in terms of the gauge (resp. support) function of the polar of $\mathcal{S}$
with \eqref{eq:gaugeASS} (resp. \eqref{eq:suppSAS}).
\begin{mytheo}
  Consider an autonomous system~\eqref{eq:autonsys}.
  The invariance of a closed convex set $\mathcal{S}$ containing the origin
  is equivalent to each of the following inequalities:
  \begin{align}
    \label{eq:suppASS}
    \forall y &\in \mathbb{R}^n, & \suppfun(\mathcal{S}, A^\Tr y) & \le \suppfun(\mathcal{S}, y).\\
    \label{eq:gaugeSAS}
    \forall x &\in \mathbb{R}^n, & \gauge(\mathcal{S}, x) & \ge \gauge(\mathcal{S}, A x).\\
    \label{eq:gaugeASS}
    \forall y &\in \mathbb{R}^n, & \gauge(\polar{\mathcal{S}}, A^\Tr y) & \le \gauge(\polar{\mathcal{S}}, y).\\
    \label{eq:suppSAS}
    \forall x &\in \mathbb{R}^n, & \suppfun(\polar{\mathcal{S}}, x) & \ge \suppfun(\polar{\mathcal{S}}, A x).
  \end{align}
  \begin{proof}
    Using \eqref{eq:ASS} and \propref{hA}, we obtain \eqref{eq:suppASS}.
    Using \eqref{eq:SAS} and \propref{gA}, we obtain \eqref{eq:gaugeSAS}.
    Using \eqref{eq:suppASS} and \propref{gaugesupp}, we obtain \eqref{eq:gaugeASS}.
    Using \eqref{eq:gaugeSAS} and \propref{gaugesupp}, we obtain \eqref{eq:suppSAS}.
  \end{proof}
\end{mytheo}

The similarity between \eqref{eq:gaugeSAS} and \eqref{eq:gaugeASS} or between \eqref{eq:suppASS} and \eqref{eq:suppSAS}
is reminescent of the fact that the stability of the system~\eqref{eq:autonsys} is equivalent to the stability of the polar system.
\begin{mycoro}
  A convex set $\mathcal{S}$ is invariant for the autonomous system~\eqref{eq:autonsys} if and only if its polar $\polar{\mathcal{S}}$
  is invariant for the following polar autonomous system:
  \(
    x_{k+1} = A^\Tr x_k.
  \)
\end{mycoro}
See also \cite{rakovic2017minkowski} for a result covering robust positively invariant sets of discrete-time linear systems with disturbances.

\subsection{Control systems}

In this section, we start by showing how to formulate the control invariance of a set
as an inequality in terms of the support function of the set or the Minkowski function of its polar.
Then we show how this allows to compute control invariant sets using a semidefine program.

The control invariance condition \eqref{eq:controlinv} can be written as
\begin{equation}
  \label{eq:controlprimal}
  \mathcal{S} \subseteq
  \begin{pmatrix}
    I_n & 0
  \end{pmatrix}
  \begin{pmatrix}
    A & B
  \end{pmatrix}^{-1}
  \mathcal{S}
\end{equation}
Note that as there is both an image an a preimage in the right-hand side.
Therefore, we cannot get an inequality in terms of $\gauge(\mathcal{S}, \cdot)$ because
of the image and cannot get an inequality in terms of $\suppfun(\mathcal{S}, \cdot)$ because of the preimage.
On the other hand, the inclusion \(
  \begin{pmatrix}
    I_n & 0
  \end{pmatrix}^{-1}
  \mathcal{S} \subseteq
  \begin{pmatrix}
    A & B
  \end{pmatrix}^{-1}
  \mathcal{S}
\)
can be rewritten as an inequality in terms of $\gauge(\mathcal{S}, \cdot)$ using \propref{gA}
but this inclusion is not equivalent to \eqref{eq:controlprimal} as we saw
in \remref{image}.

By \propref{podu}, \eqref{eq:controlprimal} is equivalent to
\begin{equation}
  \label{eq:controlpolar}
  \polar{\mathcal{S}} \supseteq
  \begin{pmatrix}
    I_n \\ 0
  \end{pmatrix}^{-1}
  \begin{pmatrix}
    A & B
  \end{pmatrix}^{\Tr}
  \polar{\mathcal{S}}
\end{equation}
See also \cite{rakovic2020robust} for a result covering robust control invariant sets of discrete-time control linear systems with disturbances.

Again,
inclusion~\eqref{eq:controlpolar} is not equivalent to
\(
  \begin{pmatrix}
    I_n & 0
  \end{pmatrix}^{\Tr}
  \polar{\mathcal{S}} \supseteq
  \begin{pmatrix}
    A & B
  \end{pmatrix}^{\Tr}
  \polar{\mathcal{S}}
\)
by \remref{image} and
there is both an image an a preimage in the right-hand side.
Hence there does not seem to be any way to use this inclusion
to write an inequality in terms of either the gauge or support function
of $\polar{\mathcal{S}}$.
We can however obtain the following theorem.
\begin{mytheo}
  \label{theo:controlprojpolar}
  The invariance of a closed convex set $\mathcal{S} \subseteq \mathcal{X}$ containing the origin
  for a control system~\eqref{eq:controlsys},
  as defined in \eqref{eq:controlinv}, is equivalent to 
  \begin{equation}
    \label{eq:controlprojpolar}
    \polar{\mathcal{S}} \supseteq
    A^\Tr
    P_{\Image(B)^{\perp}}^{\Tr}
    P_{\Image(B)^{\perp}}^{-\Tr}
    \polar{\mathcal{S}}.
  \end{equation}
  \begin{proof}
    The inclusion \eqref{eq:controlpolar} can be rewritten into:
    \begin{quote}
      For each $y \in \polar{\mathcal{S}}$,
      if $B^{\Tr}y = 0$ then $A^{\Tr}y \in \polar{\mathcal{S}}$.
    \end{quote}
    The set of $y \in \polar{\mathcal{S}}$ such that $B^{\Tr}y = 0$ is
    given by $P_{\Image(B)^{\perp}}^{\Tr} P_{\Image(B)^{\perp}}^{-\Tr} \polar{\mathcal{S}}$.
    Hence, \eqref{eq:controlpolar} is equivalent to \eqref{eq:controlprojpolar}.
  \end{proof}
\end{mytheo}

We now give an alternative approach to obtain the inclusion \eqref{eq:controlprojpolar} by
considering the algebraic system:
\begin{equation}
  \label{eq:algsys}
  Ex_{k+1} = Cx_k, x_k \in \mathcal{X}.
\end{equation}
with matrices $E, C \in \R^{r \times n}$.
We say that a set $\mathcal{S}$ is \emph{invariant} for the system~\eqref{eq:algsys} if
\begin{equation}
  \label{eq:alginv}
  C\mathcal{S} \subseteq E\mathcal{S}.
\end{equation}

The following proposition follows from an argument similar to \theoref{controlprojpolar}.
Using this argument at the level of systems instead of sets as in \theoref{controlprojpolar}
has the advantage of showing a relation between invariant sets for
two different families of systems.

\begin{myprop}[{\cite[Proposition~2]{legat2020sum}}]
  \label{prop:proju}
  The set $\mathcal{S}$ is \emph{control invariant} for the system \eqref{eq:controlsys}
  if and only if it is invariant for the system \eqref{eq:algsys}
  with $C = P_{\Image(B)^{\perp}}A$ and $E = P_{\Image(B)^{\perp}}$.
\end{myprop}

The following results follows from \propref{podu}.
Note that the inclusion~\eqref{eq:hdtahas} is exactly the inclusion~\eqref{eq:controlprojpolar}.

\begin{mytheo}[{\cite[Theorem~1]{legat2020sum}}]
  \label{theo:alginv}
  The invariance of a closed convex set $\mathcal{S} \subseteq \mathcal{X}$ containing the origin
  for an algebraic system~\eqref{eq:algsys},
  as defined in \eqref{eq:alginv}, is equivalent to 
  \begin{equation}
    \label{eq:hdtahas}
    C^{-\Tr} \polar{\mathcal{S}} \supseteq E^{-\Tr} \polar{\mathcal{S}}.
  \end{equation}
\end{mytheo}

The inclusion in \eqref{eq:alginv} can be rewritten in terms of the support
function of $\mathcal{S}$:
\begin{equation}
  \label{eq:suppalg}
  \forall y \in \R^r, \quad \suppfun(\mathcal{S}, C^\Tr y) \le \suppfun(\mathcal{S}, E^\Tr y).
\end{equation}
We can also rewrite \eqref{eq:controlprojpolar} or \eqref{eq:hdtahas} in
terms of the Minkowski function of $\polar{\mathcal{S}}$:
\begin{equation}
  \label{eq:gaugealg}
  \forall y \in \R^r, \quad \gauge(\polar{\mathcal{S}}, C^\Tr y) \le \gauge(\polar{\mathcal{S}}, E^\Tr y).
\end{equation}
The inequalities \eqref{eq:suppalg} and \eqref{eq:gaugealg} are equivalent by \propref{gaugesupp}.

The semidefinite program obtained by using an Euclidean gauge function
in \eqref{eq:gaugealg} was described in \cite{legat2018computing}
and the sum-of-squares program obtained by using the $(2d)$th root of polynomials
of degree $2d$ as gauge function was described in \cite{legat2020sum}.

We describe here the semidefinite program obtained by using piecewise
semi-Euclidean gauge functions.

\begin{myprog}
  \label{prog:piecequad}
  \begin{align}
    \notag
    \text{find}
    & \quad Q_i \succeq 0, u_{i, j} \in \R^n \quad \text{ such that}\\
      \notag
      x^\Tr C Q_i C^\Tr x & \leq x^\Tr E Q_j E^\Tr x,
      \forall x \in C^{-\Tr} \mathcal{P}_i \cap E^{-\Tr} \mathcal{P}_j,\\
      \label{eq:piecequadincl}
      & \quad \qquad \qquad \qquad \forall i, j \in \{1, \ldots, m\}\\
      \label{eq:piecequadconstr}
      v^\Tr Q_i v & \le 1, \quad \forall v \in \Vrep(\polar{\mathcal{X}}) \cap \mathcal{P}_i, \forall i \\
      \label{eq:piecequadcont}
      Q_i - Q_j & = u_{i,j} n_{ij}^\Tr + n_{ij} u_{i,j}^\Tr,
      \quad \forall i, j \in \mathcal{N}\\
      \label{eq:piecequadconv}
      Q_j n_{ij} - Q_i n_{ij} & \in \polar{(\mathcal{P}_i \cap \mathcal{P}_j)}, \quad \forall i, j \in \mathcal{N}.
  \end{align}
  where $\mathcal{N}$ and $n_{ij}$ are defined in \secref{piece}.
\end{myprog}

The constraint \eqref{eq:piecequadincl} ensures invariance,
the constraint \eqref{eq:piecequadconstr} ensures satisfaction of state-space constraints $\mathcal{X}$,
the constraint \eqref{eq:piecequadcont} ensures continuity and
the constraint \eqref{eq:piecequadconv} ensures convexity.

\begin{myrem}
  The constraint \eqref{eq:piecequadincl} can be implemented either with the
  condition provided by \propref{copos}
  or with one level of the hierarchy mentionned in \remref{copos}.
  In both cases, the resulting constraint is an LMI.
  As \eqref{eq:piecequadconstr}, \eqref{eq:piecequadcont} and \eqref{eq:piecequadconv} are linear constraints,
  \progref{piecequad} is a semidefinite program with as many LMI constraint as there are pairs $i, j$ such that
  $C^{-\Tr} \mathcal{P}_i \cap E^{-\Tr} \mathcal{P}_j$ is nonempty.
  That is, at most $m^2$ LMI constraint but there can be fewer depending on the problem data.
\end{myrem}


\begin{mytheo}
  Consider the control system~\eqref{eq:controlsys}, a
  conic partition $(\mathcal{P}_i)_{i=1}^m$ and the function
  \[
    f(x) =
      \sqrt{x^\Tr Q_i x} \quad \text{ if } x \in \mathcal{P}_i.
  \]
  The \progref{piecequad} is feasible if and only if $f$ is the
  Minkowski function of a piecewise semi-ellipsoid
  whose polar set $\mathcal{S}$
  is control invariant for the system, as defined in \eqref{eq:controlinv}.
  \begin{proof}
    By \eqref{eq:piecequadcont}, the condition \eqref{eq:piececont} holds
    and by \eqref{eq:piecequadconv}, the condition \eqref{eq:piececonv} holds.
    By \eqref{eq:piecequadconstr}, we have $\mathcal{X}^* \subseteq \mathcal{S}^*$
    hence $\mathcal{S} \subseteq \mathcal{X}$ holds.
    By \theoref{alginv} and \propref{proju}, the set $\mathcal{S}$ is control invariant for the control system.
  \end{proof}
\end{mytheo}

\begin{figure}[!ht]
  \centering
  \begin{subfigure}[b]{0.22\textwidth}
    \centering
      \includegraphics[trim=3.7cm 0cm 3.5cm 0cm, clip, width=\textwidth]{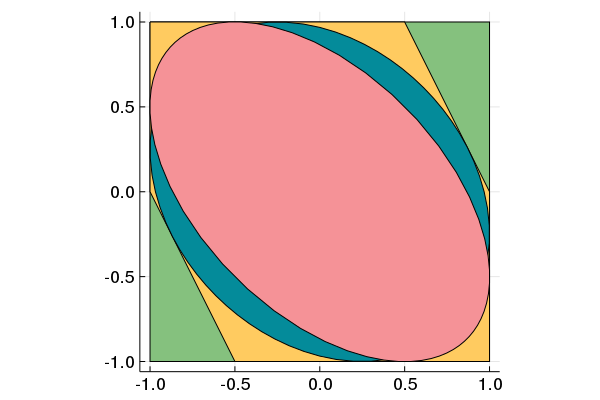}
      \caption{Primal space}
      \label{fig:SimpleControlledInvariantEllipsoid}
  \end{subfigure}%
  ~
  \begin{subfigure}[b]{0.22\textwidth}
    \centering
      \includegraphics[trim=3.7cm 0cm 3.5cm 0cm, clip, width=\textwidth]{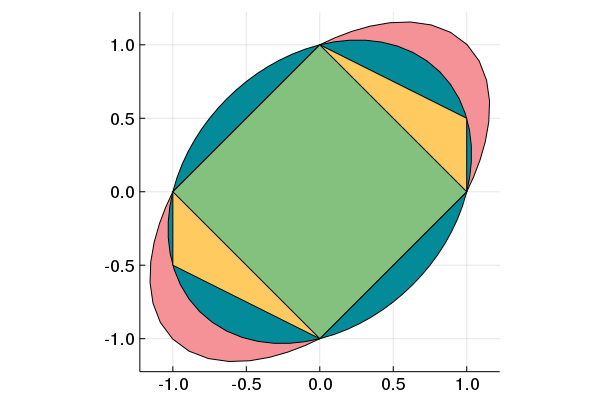}
      \caption{Polar space}
      \label{fig:SimpleControlledInvariantPolarEllipsoid}
  \end{subfigure}
  \caption{
    Figure~4 of \cite{legat2020sum} included here for completeness.
    The blue (resp. red) ellipsoid is the controlled invariant ellipsoid with
    maximal volume (resp. sum of the squares of the length of the axes).
  }
  \label{fig:simpleell}
\end{figure}

\begin{figure}[!ht]
  \centering
  \begin{subfigure}[b]{0.22\textwidth}
    \centering
      \includegraphics[trim=3.7cm 0cm 3.5cm 0cm, clip, width=\textwidth]{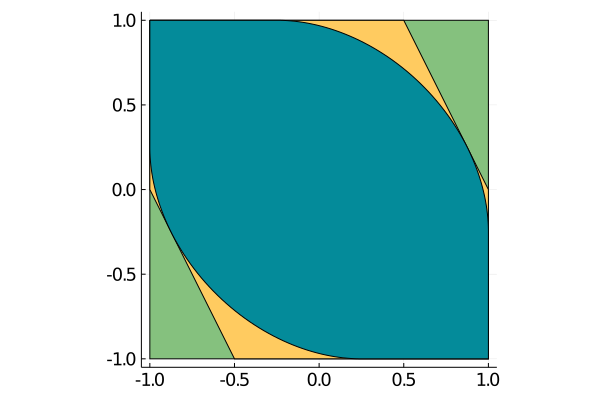}
      \caption{Primal space}
      \label{fig:SimpleControlledInvariantPiecewiseSemiEllipsoid}
  \end{subfigure}%
  ~
  \begin{subfigure}[b]{0.22\textwidth}
    \centering
      \includegraphics[trim=3.7cm 0cm 3.5cm 0cm, clip, width=\textwidth]{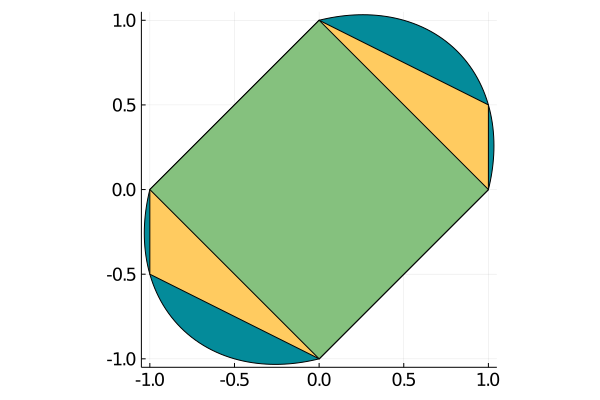}
      \caption{Polar space}
      \label{fig:SimpleControlledInvariantPolarPiecewiseSemiEllipsoid}
  \end{subfigure}
  \caption{
  Result of our algorithms applied to \exemref{simple}.
  The green set is the safe set $[-1, 1]^2$ and the blue and yellow
  sets are the optimal solutions of \progref{piecequad} with different
  conic partitions for \exemref{simple}.
%
  }
  \label{fig:simple}
\end{figure}

\begin{myexem}
  \label{exem:simple}
  Consider the simple example introduced in \cite[Example~2, 4, 5, 7]{legat2020sum}.
  The safe set is $[-1, 1]^2$ and the result of the first iteration of the standard viability kernel algorithm, 
  which is also the maximal control invariant set,
  is given by the yellow set in \figref{simple}.
  For each set, consider the conic partition given by the conic hull of each facet of the polar of the set.

  Consider the semidefinite program obtained by implementing the constraint \eqref{eq:piecequadincl} with the
  condition provided by \propref{copos} and maximizing the integral of $\suppfun^2(\mathcal{S}, x)$ over a given polytope as described in \cite{legat2020piecewiseCO}.
  This heuristic for volume maximization introduced in \cite{dabbene2017simple} is a generalization of the trace heuristic for ellipsoids. 
  The optimal solution for the first (resp. second) conic partitions is the blue (resp. yellow) set in \figref{simple}.
  Even the first solution includes both ellipsoids of \figref{simpleell}
  thanks to the ability of piecewise semi-ellipsoids
  to be polyhedral in some pieces of the spaces and ellipsoidal in other pieces.
  See \cite{legat2020piecewiseCO} for more details.
\end{myexem}

\progref{piecequad} can be adapted for a switched control system
\begin{equation*}
  x_{k+1} = A_{\sigma_k}x_k + B_{\sigma_k}u_k, x_k \in \mathcal{X}, \sigma_k \in \Sigma
\end{equation*}
with an autonomous switching signal $\sigma_k$ in a finite set $\Sigma$.
The constraint~\eqref{eq:piecequadincl} is replaced by a constraint
\[
  x^\Tr C_\sigma Q_i C_\sigma^\Tr \preceq x^\Tr E_\sigma Q_j E_\sigma^\Tr x,
  \forall i, j, x \in C_\sigma^{-\Tr} \mathcal{P}_i \cap E_\sigma^{-\Tr} \mathcal{P}_j,
\]
for each $\sigma \in \Sigma$
where $C_\sigma = P_{\Image(B_\sigma)^{\perp}}A_\sigma$ and $E_\sigma = P_{\Image(B_\sigma)^{\perp}}$.

\section{Conclusion}

We have motivated the need to consider piecewise semi-ellipsoids for
complex systems that include either state or input constraints or switched or hybrid dynamics.
This family generalizes the family of ellipsoids and polyhedra and in practice,
it allows to use the simple polyhedral representation in some part of the state space
and use an ellipsoidal surface to approximate smooth parts that would require a complex polyhedral representation.

We argued that the classical algebraic approach to control invariance does not yield any convex programming approach for computing control invariant sets of this family.
On the other hand, a geometric approach provides a convex program for control invariance by reformulating the problem in the geometric dual space.
As future work, we aim to generalize this method to sets with polynomial pieces of higher degree.
It would then be possible to enrich the approximation capabilities of the set template either
by adding more pieces or by increasing the degree of the polynomials.

In this paper, we considered the conic partition as given.
As future work, we plan to study two different approaches to obtain or refine this partition.
A first approach is to analyse the parts of the state spaces
where more halfspaces are generated after a few iterations of the standard viability kernel algorithm. This is used as a heuristic in \exemref{simple}.
A second approach is to refine the conic partition using the dual solution obtained by solving \progref{piecequad}.

\bibliographystyle{plain}
\bibliography{biblio}

\end{document}